\documentclass{amsart}

\usepackage{cite}
\usepackage{graphicx,epsfig}
\usepackage{amsmath,amsfonts,amssymb,amsthm,amsbsy}
\usepackage{color}

\theoremstyle{plain}

\theoremstyle{definition}

\begin{document}

\title[Curve shortening flow in 3D]{A short proof of monotonicity formula for curve shortening flow in 3D}

\author{Hayk Mikayelyan}

\today

\maketitle

\begin{abstract}
We give a rather elementary proof of a Huisken-type monotonicity formula
for curve shortening flow in 3D.
\end{abstract}

\thanks{Primary 35K93; Secondary 35K10 }

\thanks{Keywords: curve shortening flow, monotonicity formula}


\vspace{2cm}

We consider a closed curve in $\mathbb{R}^3$ moving by its curvature
$$
\partial_t {\bf u} = H\nu,
$$
where ${\bf u}:[0,T) \times S^1\to \mathbb{R}^3$ is the curve parametrization,
$$
H=
\frac{|{\bf u}''\times{\bf u}'|}{|{\bf u}'|^3}
$$
is the curvature and
$$\nu=
\frac{{\bf u}'\times ({\bf u}''\times{\bf u}')}{|{\bf u}'||{\bf u}''\times{\bf u}'|}
$$
is the normal vector.

Note that in this form we fix a certain parametrization which has no component
tangential to the normal $\nu$. For a general parametrization we will get
\begin{equation}\label{weaker}
\partial_t {\bf u} \cdot\nu=H
\end{equation}
and
$$
\partial_t {\bf u} \cdot\gamma=0
$$
where $\gamma=\frac{{\bf u}''\times{\bf u}'}{|{\bf u}''\times{\bf u}'|}$ is the unit binormal vector.

Assume the first singularity appears at point $0$ after finite time $T$.
We rescale the parametrization in the following way
$$
\tau=-\log(T-t), {\bf v}(\tau,x)=(T-t)^{-\frac{1}{2}}{\bf u}(t,x)
$$
and arrive at
\begin{equation}  \label{main}
\partial_\tau{\bf v}(\tau,x)
=\frac{1}{2}{\bf v}+
\frac{{\bf v}'\times ({\bf v}''\times{\bf v}')}{|{\bf v}'|^4}.
\end{equation}

\vspace{5mm}


In the paper we prove the following theorem.

{\bf Theorem.}
The solutions of (\ref{main}) satisfy the following relation
\begin{multline}\label{monformHuisk}
\frac{d}{d\tau}\int_{S^1} e^{-\frac{1}{4}|{\bf v}|^2}|{\bf v}'|dx=\\
-\frac{1}{4}\int_{S^1} | {\bf v}\cdot \gamma|^2 e^{-\frac{1}{4}|{\bf v}|^2}|{\bf v}'|dx
-\int_{S^1}|\partial_\tau {\bf v}\cdot\nu|^2 e^{-\frac{1}{4}|{\bf v}|^2}|{\bf v}'|dx,
\end{multline}
where $\nu$ is the unit normal vector and $\gamma$ the unit binormal vector.

\vspace{1cm}

This formula is the 3D codimension-2 analogue of the famous Huisken's formula
(see \cite{H2}). The proof is based on ideas introduced by Zelenjak (see \cite{Z})
for parabolic boundary value problems and
adapted by the author in 2D curve shortening context in \cite{M}.

\vspace{1cm}

{\bf Corollary 1.}
Stationary solutions of (\ref{main}) are plane curves.

\vspace{1cm}

{\bf Corollary 2.}
Let us set
$$
E(\tau)=\int_{S^1} e^{-\frac{1}{4}|{\bf v}|^2}|{\bf v}'|dx\,\,\,\,\,\,\text{and}\,\,\,\,\,\,
\Pi(\tau)=\frac{1}{4}\int_{S^1} | {\bf v}\cdot \gamma|^2 e^{-\frac{1}{4}|{\bf v}|^2}|{\bf v}'|dx.
$$
Then
$$
\Pi(\tau)\leq -E'(\tau).
$$

\vspace{1cm}

{\bf Proof of the Theorem.}

For the system (\ref{main}) let us try to obtain a monotonicity formula of the form
\begin{equation}\label{monform}
\frac{d}{d\tau}\int_{S^1} F(v_1,v_2,v_3,v_1',v_2',v_3')dx\leq
-\int_{S^1}|\partial_\tau {\bf v}\cdot\nu|^2\rho(v_1,v_2,v_3,v_1',v_2',v_3')dx,
\end{equation}
where $\rho$ is positive.

Differentiating the left hand side of (\ref{monform}) and integrating by parts
we get
\begin{align}\label{tau1}
\partial_\tau v_1\left[\frac{\partial F}{\partial\xi_1}-
\frac{\partial^2 F}{\partial\xi_1\partial\eta_1}v_1'-
\frac{\partial^2 F}{\partial\xi_2\partial\eta_1}v_2'-
\frac{\partial^2 F}{\partial\xi_3\partial\eta_1}v_3'{\color{blue}-
\frac{\partial^2 F}{\partial\eta_1^2}v_1''-
\frac{\partial^2 F}{\partial\eta_1\partial\eta_2}v_2''-
\frac{\partial^2 F}{\partial\eta_1\partial\eta_3}v_3''}
\right]+\\
\label{tau2}
\partial_\tau v_2\left[\frac{\partial F}{\partial\xi_2}-
\frac{\partial^2 F}{\partial\xi_1\partial\eta_2}v_1'-
\frac{\partial^2 F}{\partial\xi_2\partial\eta_2}v_2'-
\frac{\partial^2 F}{\partial\xi_3\partial\eta_2}v_3'{\color{blue}-
\frac{\partial^2 F}{\partial\eta_1\partial\eta_2}v_1''-
\frac{\partial^2 F}{\partial\eta_2^2}v_2''-
\frac{\partial^2 F}{\partial\eta_2\partial\eta_3}v_3''}
\right]+\\
\label{tau3}
\partial_\tau v_3\left[\frac{\partial F}{\partial\xi_3}-
\frac{\partial^2 F}{\partial\xi_1\partial\eta_3}v_1'-
\frac{\partial^2 F}{\partial\xi_2\partial\eta_3}v_2'-
\frac{\partial^2 F}{\partial\xi_3\partial\eta_3}v_3'{\color{blue}-
\frac{\partial^2 F}{\partial\eta_1\partial\eta_3}v_1''-
\frac{\partial^2 F}{\partial\eta_2\partial\eta_3}v_2''-
\frac{\partial^2 F}{\partial\eta_3^2}v_3''}
\right].
\end{align}

In the right hand side of (\ref{monform}) using (\ref{main}) we obtain
\begin{align} |\partial_\tau{\bf v}\cdot\nu|^2=
\left[\left(
\begin{array}{cc} \partial_\tau v_1\\\partial_\tau v_2\\\partial_\tau v_3
\end{array}
\right)\cdot \nu \right]\left[\left(  \frac{1}{2}{\bf v}+
\frac{{\bf v}'\times ({\bf v}''\times{\bf v}')}{|{\bf v}'|^4} \right)\cdot\nu\right]=
\\
\left[\left(
\begin{array}{cc} \partial_\tau v_1\\\partial_\tau v_2\\\partial_\tau v_3
\end{array}
\right)\cdot \nu \right]\left[  \frac{1}{2}{\bf v}\cdot\nu\right]
+
\left[\left(
\begin{array}{cc} \partial_\tau v_1\\\partial_\tau v_2\\\partial_\tau v_3
\end{array}
\right)\cdot \nu \right]\left[
\frac{{\bf v}'\times ({\bf v}''\times{\bf v}')}{|{\bf v}'|^4} \cdot\nu\right]=
\\\label{gammasecondterm}
\left[\left(
\begin{array}{cc} \partial_\tau v_1\\\partial_\tau v_2\\\partial_\tau v_3
\end{array}
\right)\cdot \nu \right]\left[  \frac{1}{2}{\bf v}\cdot\nu\right]
+{\color{blue}
\left[\left(
\begin{array}{cc} \partial_\tau v_1\\\partial_\tau v_2\\\partial_\tau v_3
\end{array}
\right)\cdot
\frac{{\bf v}'\times ({\bf v}''\times{\bf v}')}{|{\bf v}'|^4}\right]}.
\end{align}
Observe that
$$
\frac{{\bf v}'\times ({\bf v}''\times{\bf v}')}{|{\bf v}'|^4}=
\frac{1}{|{\bf v}'|^4}{\color{blue}
\left(
\begin{array}{cc}
v_1''({v'}_2^2+{v'}_3^2)-v_2''v_1'v_2'-v_3''v_1'v_3'
\\
-v_1''v_1'v_2'+v_2''({v'}_1^2+{v'}_3^2)-v_3''v_2'v_3'
\\
-v_1''v_1'v_3'-v_2''v_2'v_3'+v_3''({v'}_2^2+{v'}_3^2)
\end{array}
\right)}
$$
We require that
\begin{equation}\label{Frho}
D^2_\eta F(\xi, \eta)=\frac{\rho(\xi,\eta)}{|\eta|^4}\left(
\begin{array}{ccc}
(\eta_2^2+\eta_3^2) & -\eta_1\eta_2 & -\eta_1\eta_3
\\
-\eta_1\eta_2 & (\eta_1^2+\eta_3^2) & -\eta_2\eta_3
\\
-\eta_1\eta_3 & -\eta_2\eta_3 & (\eta_1^2+\eta_2^2)
\end{array}
\right)=\rho(\xi,\eta)|\eta|^{-1}D^2|\eta|,
\end{equation}
which means that the last three terms in (\ref{tau1}), (\ref{tau2}) and (\ref{tau3}) will ``take care'' of the second term in (\ref{gammasecondterm}).

Setting $\nu=(\nu_1,\nu_2,\nu_3)^T$ we further require that the first term in (\ref{gammasecondterm}) multiplied by $-\rho(\xi,\eta)$ be larger or equal than
to the remaining terms in (\ref{tau1}), (\ref{tau2}) and (\ref{tau3}).

\begin{align}\label{geqFrho}
-\rho(\xi,\eta)
\left[\left(
\begin{array}{cc} \partial_\tau v_1\\\partial_\tau v_2\\\partial_\tau v_3
\end{array}
\right)\cdot \nu \right]\left[  \frac{1}{2}{\bf v}\cdot\nu\right]=
-\frac{\rho}{2}\partial_\tau{\bf v}\cdot[( {\bf v}\cdot \nu ) \nu]\geq
\\
\partial_\tau v_1\left[\frac{\partial F}{\partial\xi_1}-
\frac{\partial^2 F}{\partial\xi_1\partial\eta_1}v_1'-
\frac{\partial^2 F}{\partial\xi_2\partial\eta_1}v_2'-
\frac{\partial^2 F}{\partial\xi_3\partial\eta_1}v_3'
\right]+\\
\partial_\tau v_2\left[\frac{\partial F}{\partial\xi_2}-
\frac{\partial^2 F}{\partial\xi_1\partial\eta_2}v_1'-
\frac{\partial^2 F}{\partial\xi_2\partial\eta_2}v_2'-
\frac{\partial^2 F}{\partial\xi_3\partial\eta_2}v_3'
\right]+\\
\partial_\tau v_3\left[\frac{\partial F}{\partial\xi_3}-
\frac{\partial^2 F}{\partial\xi_1\partial\eta_3}v_1'-
\frac{\partial^2 F}{\partial\xi_2\partial\eta_3}v_2'-
\frac{\partial^2 F}{\partial\xi_3\partial\eta_3}v_3'
\right].
\end{align}
\vspace{5mm}

Now let us observe that the functions
$$
F(\xi,\eta)=\rho(\xi,\eta)=|\eta|e^{-\frac{|\xi|^2}{4}}
$$
satisfy (\ref{Frho}) and inequality (\ref{geqFrho}). (\ref{Frho})
is obvious and to see that
 (\ref{geqFrho}) is satisfied we compute
 \begin{align}
\left(
\begin{array}{cc}
\frac{\partial F}{\partial\xi_1}-
\frac{\partial^2 F}{\partial\xi_1\partial\eta_1}\eta_1-
\frac{\partial^2 F}{\partial\xi_2\partial\eta_1}\eta_2-
\frac{\partial^2 F}{\partial\xi_3\partial\eta_1}\eta_3
\\
\frac{\partial F}{\partial\xi_2}-
\frac{\partial^2 F}{\partial\xi_1\partial\eta_2}\eta_1-
\frac{\partial^2 F}{\partial\xi_2\partial\eta_2}\eta_2-
\frac{\partial^2 F}{\partial\xi_3\partial\eta_2}\eta_3
\\
\frac{\partial F}{\partial\xi_3}-
\frac{\partial^2 F}{\partial\xi_1\partial\eta_3}\eta_1-
\frac{\partial^2 F}{\partial\xi_2\partial\eta_3}\eta_2-
\frac{\partial^2 F}{\partial\xi_3\partial\eta_3}\eta_3
\end{array}
\right)=\\
-\frac{1}{2}e^{\frac{|\xi|^2}{4}}\left[ |\eta|\xi -\frac{\xi\cdot \eta}{|\eta|}\eta\right]=
-\frac{\rho}{2}\left[ \xi -\frac{\xi\cdot \eta}{|\eta|^2}\eta\right]
\end{align}
We need to check that
\begin{align}
-\frac{\rho}{2}\partial_\tau{\bf v}\cdot[( {\bf v}\cdot \nu ) \nu]=
-\frac{\rho}{2}\left[ \frac{1}{2}\xi+H\nu\right]\cdot[( \xi\cdot \nu ) \nu])=
\\
-\frac{\rho}{4}| \xi\cdot \nu |^2-\frac{H\rho}{2} \xi\cdot \nu \geq
\\
-\frac{\rho}{2}\left[ \frac{1}{2}\xi+H\nu\right]\cdot \left[ \xi -\frac{ \xi\cdot \eta }{|\eta|^2}\eta\right]=
-\frac{\rho}{4}|\xi|^2
+\frac{\rho}{4}\frac{|\xi\cdot \eta|^2}{|\eta|^2}-\frac{H\rho}{2}\xi\cdot \nu.
\end{align}
The result follows from mutual orthogonality of $\nu$, $\eta$ and $\gamma$, and the equality
$$
|\xi|^2=| \xi\cdot \nu |^2+\frac{| \xi\cdot \eta |^2}{|\eta|^2}+|\xi\cdot \gamma |^2.
$$

\vspace{1cm}

{\bf Remark.}
In the proof we use the compactness of the curve only when doing partial integration
to obtain (\ref{tau1})-(\ref{tau3}). For a non-compact curve parametrized by $x\in(-\infty,\infty)$
all the calculations will work provided the integrals are finite and
$$
\lim_{x\to\pm\infty} \partial_{\eta_j}F({\bf v},{\bf v}')\partial_\tau v_j=0
$$
for $j=1,2,3$.










\vspace{5mm}

\address{\noindent Mathematcal Sciences\\
The University of Nottingham Ningbo China\\
Taikang Dong Lu Nr. 199,
Ningbo 315100\\ PR China}
\vspace{2mm}

\email{\noindent Hayk.Mikayelyan@nottingham.edu.cn}

\end{document}